# ON DIFFUSION APPROXIMATION OF A SLOW COMPONENT FOR SOLUTION OF STOCHASTIC DIFFERENTIAL EQUATION OF ITO: THEORY AND ITS APPLICATION TO THE GROUNDING OF THE FICK'S DIFFUSION LAWS

## Doobko V.A.

Consider the system of stochastic differential equations:

$$\varepsilon dy_\varepsilon(t) = -A(x_\varepsilon(t);t)y_{\varepsilon(t)}dt + F(x_\varepsilon(t);t)dt + B(x_\varepsilon(t);t)dw(t),$$
$$dx_\varepsilon(t) = C(x_\varepsilon(t);t)y_\varepsilon(t)dt + Q(x_\varepsilon(t);t)dt + P(x_\varepsilon(t);t)dw(t),\quad(1)$$

where $w(t)$ - m-dimensional Winner's process, $y, F \in R^n; x, Q \in R^{n'}$, $A, B, C, P$ - are matrices of size $n \times n, n \times m, n' \times n, n' \times m$ respectively, and $\varepsilon$ - is a small parameter.

That type of equations like (1) we discussed in [1], for the case of $A = A(x; y; t), B = B(x; y; t)$ and limit equations for the slow components are found.

The Brownian motion is an example of the process described by (1)

$$\varepsilon dv_\varepsilon(t) = -A(x_\varepsilon(t);t)v_\varepsilon(t)dt + F(x_\varepsilon(t);t)dt + B(x_\varepsilon(t);t)dw(t),$$
$$dx_\varepsilon(t) = v_\varepsilon(t)dt \quad (2)$$

where $x(t) \in R^3$, $v(t) \in R^3$ are coordinates and velocity of a Brownian particle. In this work we study system (2) in more details.

It is known that a diffusion motion can be described as a Markov process in $R^3(x)$:

$$dx(t) = q_0(x(t);t)dt + C_0(x(t);t)dw(t) \quad (3)$$

This equation corresponds to the following one for density $\rho(x;t)$

$$\frac{\partial}{\partial t}\rho(x;t) = \frac{\partial}{\partial x_i}(q_{0,i}(x;t)\rho(x;t)) + \frac{\partial^2}{\partial x_i \partial x_j}(C_{0,i,k}(x;t)C_{0,j,k}(x;t)\rho(x;t)) \quad (4)$$

Here and below the summation over the repeating indices is assumed.

**Theorem** (see [1]). If in equations (1):
$Q(x;t), F(x;t), B(x;t)$ - are continuous for the variable $t$ and bounded and satisfy Lipchitz's ($\forall x \in R^n$) condition.

1. $C_{i,e}(x;t), a_{i,j}(x;t)$ - are bounded together with their first and second derivatives in $x, t$.
2. There exists $\gamma > 0$, such that for $(Z, AZ) \geq \gamma |Z|^2, Z \in R^n$.
3. $\varepsilon^3 M|y(0)|^4 = 0(\varepsilon), x_\varepsilon(0) = x(0), y_\varepsilon(0) = y(0),$

- then $\exists l.i.m._{\varepsilon \to 0} x_\varepsilon(t) = x(t)$, where $x(t)$ is the solution of the approximation equation (3) with the following coefficients:

$$q_{0,i}(x;t) = (C(x;t)A^{-1}(x;t)F(x;t))_i + q_i(x;t) + [p_{i,k}(x;t)b_{l,k}(x;t) \times$$

$$\frac{\partial}{\partial x_j}(C(x;t)A^{-1}(x;t))_{i,l}] + C_{j,l}(x;t)\left[\frac{\partial}{\partial x_j}(C(x;t)A^{-1}(x;t))_{i,l'}\right] \times$$

$$\left[\int_0^\infty (\exp\{-\lambda A(x;t)\}B_k(x;t)_l \times (\exp\{-\lambda A(x;t)\}B_k(x;t)_{l'}d\lambda\right];$$

$$C_0(x;t) = C(x;t)A^{-1}(x;t)B(x;t) + P(x;t),$$

where $B_k$ - vector, composed from the matrix's line elements $B(x;t)$.

The requirements to matrices and their coefficients in the Theorem give us only uniqueness of (1), (2) and allow rewriting (2) as follows (see condition 3)

$$x_\varepsilon(T) - x(0) = -\varepsilon \int_0^T C(x_\varepsilon(t);t) A^{-1}(x_\varepsilon(t);t) dy_\varepsilon(t) + \int_0^T [C(t,x_\varepsilon(t)) A^{-1}(x_\varepsilon(t);t) F(x_\varepsilon(t);t) +$$
$$+ Q(x_\varepsilon(t);t)] dt + \int_0^T [C(x_\varepsilon(t);t) A^{-1}(x_\varepsilon(t);t) B(x_\varepsilon(t);t) + P(x_\varepsilon(t);t)] dw(t)$$

The proof of the theorem is based on the following scheme:

1) We establish that at conditions of the Theorem, $\forall t \geq 0$ $\varepsilon^3 M|y_\varepsilon(t)|^4 = O(\varepsilon)$

2) We integrate part by part the right component of the equation (4), using Ito formula and then the transformations are carried with the help of the equation for the integral from $y_{\varepsilon,i}(t) \cdot y_{\varepsilon,j}(t)$.

3) We evaluate the proximity of integral equation solution, obtained from equation (4) after operation 2) and off of all components putt estimated to be approximately equal to $\varepsilon^\delta$, $\delta > 0$ with the solution $x(t)$ (in a sense of l.i.m.).

Return to (2) for Brown motion but in more simple form:
$$mdv(t) = -6\pi\mu(x(t))rv(t)dt + b(x(t))dw(t)$$
$$dx(t) = v(t)dt \qquad (5a)$$

where - $x(t)$ - coordinate, $v(t)$ - velocity, $m$ – mass and $r$ – radius of Brownian particles, $\mu(x)$ - viscosity, $w(t)$ - Winner's process $\in R^3$; $\mu(x)$ and $b(x;t)$ - scalar functional.

With the help of this example, we will examine algorithm 1) – 3), and estimate time interval permitting to yield the Fick's diffusion laws with the help of the Theorem.

Assume that temperature of the medium is stable. It is equivalent to the following:
$$\lim M[d_t|v^2(t)|^2 / x(t_1) = x] = O \cdot dt$$
$$t \to \infty, \quad (t_1 < t)$$

Using equation (5) we obtain the following formula:
$$\lim_{t \to \infty} M\left[\frac{|v(t)|^2 m}{2} / x(t) = x\right] = \frac{3}{4} \frac{b^2(x)}{6\pi\mu(x)} = \frac{2}{3} K_B T^0 \qquad (6)$$

where $K_B$ - Boltzmann constant, $T^0$ - temperature of the medium (Kelvin scale)

This formula allows to represent equation (5a) as follows
$$mdv(t) = -6\pi\mu(x(t))rv(t)dt + (12\pi K_B T^0 \mu(x) r)^{1/2} dw(t) \qquad (5b)$$

Let's examine dimensionless variables of the (5b) equation
$$m_0 = \frac{m}{6\pi\bar\mu\ell\tau}; \quad s = \frac{m_0 t}{\tau}; \quad u(s) = \frac{\tau v(t)}{m_0 \ell}; \quad y(s) = \frac{x(t)}{\ell}; \quad \eta(y) = \frac{\mu(y\ell)}{\bar\mu} \leq 1,$$

where $\bar\mu = \sup_x \mu(x)$, $\tau$ - is characteristic time interval between successive observation; $\ell$ - selected dimensional scale.

The equation (5b) with the new variables is substituted by
$$\varepsilon \, du(s) = -\eta(y(s)) u(s) ds + \sigma(y(s)) dw(s) \qquad (7a)$$
$$dy(s) = u(s) ds \qquad (7b)$$

where
$$\sigma(y) = [2m^{-1} K_B T^0 \tau^2 \ell^{-1} r^{-1} \eta(y)]^{1/2}$$
$$\varepsilon = m_0^2 \cdot \ell \cdot r^{-1} \qquad (7c)$$

The next task is investigation of the Theorem conclusions for its applicability to the actual Brownian process. For that purpose we should introduce particular physical values for corresponding components in equations (5b), (7a-7c).

If $\ell = 1 meter$, $\tau = 1 sec.$, $\bar{\rho}$ (density of particle) $= \bar{\rho}$ (density of medium)$=10^3 kg/\ell^3$ (water density), $r = 10^{-6}\ell$, $T^0 = 300^0 K$, $v = \bar{\mu}/\bar{\rho} = 1,4 \cdot 10^{-6} \ell^2 sec^{-1}$ (kinematics' viscosity), we obtain the following values:

$$m_0 = m(6\pi\mu\ell\tau)^{-1} \cong 4r^3(3v\ell\tau)^{-1} \cong 10^{-12};$$
$$\bar{\sigma}^2 = \sup \sigma^2(y)\sigma^2(y) \cong 10; \varepsilon \cong 10^{-18} \quad (8)$$

I.e. $\varepsilon$ - may really be considered as the small parameter for (7a) equation.

These numerical values can be used for estimation of the moments :

$$\varepsilon^n \lim_{S \to \infty} M|u(s)|^{2n} = \varepsilon^n |u(\infty)|^{2n}, n = 1, 2, 3.$$

From the equation (7a) we have:

$$\varepsilon^n M|u(s)|^{2n} = (\int_0^S \exp\left\{-2n\varepsilon^{-1}\int_{s_1}^s M[\eta(y(s_2))|u(s_2)|^{2n}](M|u(s_2)|^{2n})^{-1} ds_2\right\} \times$$
$$\times \varepsilon^{(n-2)} n(2n+1) M[\sigma^2(y(s_1))|u(s_1)|^{2(n-1)}]ds_1) + \quad (9)$$
$$+ \varepsilon^n M|u(0)|^{2n} \exp\left\{-2n\varepsilon^{-1}\int_0^s M[\eta(y(s_1)|u(s_1)|^{2n}](M|u(s_1)|^{2n-1})^{-1}\right\}$$

Having defined values for $n = 1,2,3$ we obtain:

$$\varepsilon M|u(\infty)|^2 \le \sup_y \left(3\sigma^2(y)2^{-1}\eta(y)\right) \le 2^{-1}30 = K_1, \quad (10a)$$
$$\varepsilon^2 M|u(\infty)|^4 \le K_1^2, \quad \varepsilon^3 M|u(\infty)|^6 \le K_1^3,$$

- due to the condition (6) and equation (9).

If we suppose, that Brownian particles system at the initial time was in a state of thermodynamic equilibrium with the medium, these values should be observed at $t = 0$. Thus we assume that

$$\varepsilon^2 M|u(0)|^6 \le K_1^3 \quad (10b).$$

Let us suppose that $\mu(x)$ is a limited function with its secondary derivatives on $x, \forall x \in R^3$. In this case, due to the Theorem's conclusion we can affirm the existence of $[0, s)$ interval, in which the solution of $y(s)$ equation

$$d\bar{y}(s) = -2^{-1}\left[\sigma(\bar{y}(s)\eta^{-1}(\bar{y}(s))\right]^2 \nabla_{\bar{y}} \ln \eta(\bar{y}(s))ds + \sigma(\bar{y}(s))\eta^{-1}(\bar{y}(s))dw(s),$$
$$\bar{y}(0) = y(0) \quad (11)$$

- is a diffusion approximation for $y(s)$ process.

Let us examine the closeness between the solutions $y(s)$ and $\bar{y}(s)$ of equations (7b) and (11) with the help of our algorithm. We obtain the following equality:

$$M[y(s) - \bar{y}(s)]^2 = M[-\varepsilon u(s)\eta^{-1}(y(s_1))\big|_0^s - 2^{-2}\varepsilon^{-2}u(s_1)(u(s_1), \nabla_y \eta^{-2}(y(s_1))\big|_0^s +$$
$$2^{-2}\varepsilon^2 \int_0^s u(s_1)(u(s_1), \nabla_y)^2 \eta^{-2}(y(s_1))ds_1 + 2^{-2}\varepsilon \int_0^s \{\sigma(y(s_1))dw(s_1)(u(s_1), \nabla_y)\eta^{-2}(y(s_1)) +$$
$$\sigma(y(s_1))u(s_1)(dw(s_1), \nabla_y \eta^{-2}(y(s_1))\} + 2^{-2} \int_0^s \{\sigma^2(y(s_1))\nabla_y \eta^{-2}(y(s_1)) - \sigma^2(\bar{y}(s_1))\nabla_{\bar{y}} \eta^{-2}(\bar{y}(s_1))\}ds_1 +$$
$$\int_0^s [\sigma^2(y(s_1))\eta^{-2}(y(s_1)) - \sigma(\bar{y}(s_1))\eta^{-1}(\bar{y}(s_1))]dw(s_1)]^2.$$

Basing on the values for moment (9), (10) and condition (8), we can make sure that $\forall s \in [0;1)$:

$$M[y(s)-\bar{y}(s)]^2 \le \varepsilon K_2(s) + \int_0^s K_3(s)M[y(s)-\bar{y}(s_1)]^2 ds_1,$$

where

$$K_2(s) = 2K_1[2\sup_y|\eta(y)|^{-2} + \varepsilon K_1\sup_y|\nabla_y\eta^{-2}(y)|^2 + \varepsilon s^2 K_1\sup_y|\nabla_y\eta^{-2}(y)|^2 + \tag{12}$$
$$+ \varepsilon s^2 K_1\sup_y|\nabla_y\eta^{-2}(y)|^2 + s\sup_y|\sigma(y)|^2|\nabla_y\eta^{-2}(y)|^2],$$

$$K_3(s) = 3\max\{\sup_{y,\Delta} s(\sigma^2(y)\nabla_y\eta^{-2}(y) - \sigma^2(y+\Delta)\nabla_y\eta^{-2}(y+\Delta))^2\Delta^{-2});$$
$$\sup_{y,\Delta}(|\sigma(y)\eta^{-1}(y) - \sigma(y+\Delta)\eta^{-1}(y+\Delta)|^2\Delta^{-2})\}. \tag{13}$$

According to the Gronuall-Bellman inequality $\forall s \in [0;1)$

$$M|y(s)-\bar{y}(s)|^2 \le \varepsilon \bar{K}_2 \exp\{\bar{K}_3 s\},$$
$$\bar{K}_2 = K_2(s=1), \quad \bar{K}_3 = K_3(s=1),$$

or in initial variables $|\Delta_0|^2 = M|x(t)-\bar{x}(t)|^2 \le \ell^2 \bar{K}_2 10^{-18} \exp\{\bar{K}_3 \cdot \tau^{-1} \cdot t \cdot 10^{-12}\}$.

It is reasonable to compare solution for stochastic processes in case if $\bar{x}(t)$ is scattering on coordinate variable due to the random component in equation (11) which exceed the error of approximation in the range of supervision, i.e.

$$|\Delta\bar{x}|^2 = M\left|\bar{x}(t) - \int_o^t g(\bar{x}(t_1))dt_1 - x(0)\right|^2 > \delta^2, \tag{14}$$

$$\delta^2 > |\Delta_0|^2, \tag{15}$$

where $g(x) = 2^{-2}r^{-1}b^2(x)\nabla_x\mu^{-2}(x)$ - shift vector.

Brownian particle's position as a single unit may be determined with the precision not exceeding its diameter. Therefore, let us select $\delta^2 = r^2 = 10^{-12}\ell^2$, as a minimum possible error of measurement.

As appears from the equation (11),

$$M\left|\bar{y}(s) - u(0) - 2^{-2}\int_0^s \sigma^2(\bar{y}(s_1))\nabla_{\bar{y}}\eta^{-2}(\bar{y}(s_1))ds_1\right| \ge K_4 s,$$

where $K_4 = 3\inf_y|\sigma^2(y)\eta^{-2}(y)|$, or in initial variables

$$M\left|\bar{x}(t) - x(0) - \int_0^t g(\bar{x}(t_1))dt_1\right|^2 \ge \ell^2 t 10^{-12} K_4 \tau^{-1}.$$

The condition (14) is observed as soon as $\tau^{-1}t\bar{K}_3 > 1$. Condition (15) is observed if $\forall s \in [0;1)$ $1 > \bar{K}_2 10^{-6}\exp\{1\}$ and (14) is equivalent to the statement that corresponding values in time interval $t \in \tau[\bar{K}_4^{-1} = 300^{-1}, 10^{12}\bar{K}_3^{-1})$ defining $K_2, K_3$ (equations (12), (13)) may be approximately equal to $10^6$.

Assuming that $K_2 = 10^2$ which conforms to initial values in (8) for the model under study, we come to the conclusion that in this case time interval between accurate and approximate solutions, $x(s)$, is not less than $10^{10}$ sec.

Kolmogorov's equation for density corresponding to (4) is:

$$\frac{\partial}{\partial t}\rho(x;t) = \nabla_x[(D(x)\nabla_x \ln \eta(x))\rho(x;t) + \nabla_x D(x)\rho(x;t)]$$

where:

$$D(x) = K_B T^0 (6\pi\mu(x)r)^{-1} \qquad (T^0 = \text{const}),$$

and it may be transformed as follows [1]:

$$\frac{\partial}{\partial t}\rho(x;t) = \nabla_x \left[ D(x)\nabla_x \rho(x;t) \right].$$

This equation corresponds to mathematical formulation of phenomenological diffusion law for inhomogeneous medium, which conforms well to the second Fick's law. Therefore, model (5) seems to be correct in reflecting the actual dynamics of a Brownian particle. The knowledge of these estimations gives the possibility to judge about adequacy of the initial assumption about dynamics of the real processes.

Note, that some type of equations (1) were discussed in [1], for the case of $A = A(x;y;t)$, $B = B(x;y;t)$ and also limiting theorems are formulated for the slow component $x_\varepsilon(t)$.

## REFERENCES

1. Дубко В.А. Метод диффузионной аппроксимации в исследовании моделей стохастических динамических систем. – Владивосток: Дальнаука, 1994. 107 с.
( Doobko V.A. Method of diffusion approximation in the study of models of stochastic dynamic systems. – Vladivostok: Dalnauka, 1994, 107 p. (in Russian))

**Abstract.** For the concrete model of Brownian particles dynamics in non-uniform environment, the time interval estimation is constructed, on which phenomenological Fick laws for diffusion phenomenon description can be used. The knowledge of these estimations gives the possibility to judge about adequacy of the initial assumption about dynamics of the real processes. Noted, that such an assessment has not been introduced yet, up to date.